\documentclass[12pt,a4paper]{article}
\usepackage[utf8]{inputenc}
\usepackage{amsfonts,amsthm,amssymb,amsmath}

\newcommand {\Q} {\mathbb{Q}}
\newcommand {\R} {\mathbb{R}}

\newcommand {\Z} {\mathbb{Z}}
\newcommand {\FF} {\mathfrak{F}}
\newcommand {\II} {\mathfrak{I}}
\newcommand {\QQ} {\mathbb{Q}_*}
\newcommand {\QF} {\mathbb{Q}_{*fin}}
\newtheorem{theor}{Theorem}
\newtheorem{lemma}{Lemma}

\title {\flushright{\small{MSC 03H05  }}\\ 
  \center{On The Model Of Hyperrational Numbers With Selective Ultrafilter}}
\author {A.~Grigoryants \\ Moscow State University Yerevan Branch}
\begin {document}

\maketitle
\abstract{
  In standard construction of hyperrational numbers using ultrapower we assume that the ultrafilter is selective.
  It makes possible to assign real value to any finite hyperrational number. So, we can consider hyperrational numbers with 
  selective ultrafilter as extension of traditional real numbers. 
  Also proved the existence of strictly monotonic or stationary representing sequence for any hyperrational number.
}


\section {Notation and definitions}
We use standard set theoretic notation (see \cite{jech}). Let us give some well-known definition for convenience.

\emph{Partition} of a set $S$ is a pairwise disjoint family $\{S_i\}_{i\in I}$ of nonempty 
subsets such that $\bigcup_{i\in I} S_i = S$.

Let $\FF\subset 2^\omega$ be non-principal ultrafilter on $\omega$. We'll call elements in $\FF$ \emph{big subsets} 
(relative to $\FF$) and elements not in $\FF$ \emph{small subsets} (relative to $\FF$).

$\FF$ is called \emph{selective ultrafilter} if for every
partition $\{S_n\}_{n\in\omega}$ of $\omega$ into $\aleph_0$ pieces such that $S_n\notin\FF$ for all $n$ (\emph{small partition})
there exists $B\in\FF$ (\emph{big selection}) such that $B\cap S_n$  is singleton for all $n\in\omega$.
Equivalently, $\FF$ is selective if for every function $f:\omega\to\omega$ such that $f^{-1}(i)$ is small for every
$i$, there exists a big subset $B$ such that restriction $f|_B$ is injective.

The restriction of $\FF$ to a big subset $B\subset\omega$ defined as $\FF_B = \{J\cap B~|~J\in\FF\}\subset 2^B$
is selective ultrafilter on $B$.

Continuum hypothesis (CH) implies existence of selective ultrafilters due to Galvin \cite[theor. 7.8]{jech}.
The result of Shelah \cite{shelah} shows that existence of selective ultrafilters is unprovable in ZFC.
So, we continue with ZFC \& CH to ensure the existence of selective ultrafilter.

Let $[S]^k = \{X\subset S:|X|=k\}$ is the set of all subsets of $S$ that have exactly $k$ elements.
If $\{X_i\}_{i\in I}$ is a partition of $[S]^k$ then a subset $H\subset S$ is \emph{homogeneous} for the
partition if for some $i:~[H]^k\subset X_i$.

The following fact is special case of Kunen's theorem proven in \cite[theor. 9.6]{com_neg}:
\begin{theor}\label{t0}
  An ultrafilter on $\omega$ is selective if and only if for every partition of $[\omega]^2$ into two pieces there is
  a big homogeneous set.
\end{theor}

Due to this theorem selective ultrafilters are also called Ramsey ultrafilters.
We'll give simplified proof of theorem \ref{t0}.

We say that filter $\FF\subset 2^\omega$ is \emph{normal} if for any collection $\{A_i\}_{i\in\omega}\subset\FF$
there exists $B\in\FF$ such that for any $i,j\in B: i<j\implies j\in A_i$. Equally, we can say that $\FF$ is normal if
there exists $I\in\FF$ such that above definition holds for any collection $\{A_i\}_{i\in I}\subset\FF$.
Indeed we can expand given collection adding $A_k=\omega$ for $k\notin I$, then apply definition to $\{A_i\}_{i\in\omega}$
and intersect obtained $B$ with $I$.

We continue with fixed selective ultrafilter $\FF$ on $\omega$.

Let us denote $\QQ = \Q^\omega/\sim_\FF$ the quotient of $\Q^\omega = \{(x_i)_{i\in\omega}~|~\forall i: x_i\in\Q \}$
by the following equivalence relation:
$$(x_i) \sim_\FF (y_i) \Leftrightarrow \{i\in I ~|~ x_i = y_i\} \in\FF ~.$$
This is well-known ultrapower construction widely used in model theory 
and in Robinson's non-standard analysis (see \cite{rob0}, \cite{rob1}).
We only added the property of selectivity to $\FF$. So, we call elements of $\QQ$ \emph{hyperrational numbers}. 
There is natural embedding $\iota: \Q\to\QQ$ where $\iota(q)$ is the equivalence class of constant sequence
$(q,q,\dots,q,\dots)$. We'll identify $\Q$ with $\iota(\Q)$.
Also $\QQ$ satisfies the transfer principle. So, all true first order statements about $\Q$ are also valid in $\QQ$.
In particular, $\QQ$ is ordered field.

We call element $x\in\QQ$ \emph{infinitely large} if $|x|>|q|$ for all $q\in\Q$, 
\emph{infinitesimal} if $|x|<|q|$ for all $q\in\Q$. Otherwise $x$ is called \emph{finite}.
We write $x<\infty$ if $x$ is finite or infinitesimal.

We use short notation $x_n$ for some representative  $(x_n)_{n\in\omega}$ of equivalence class $x = [(x_n)_{n\in\omega}]\in\QQ$.
So, we can make arbitrary changes to sequence $x_n$ on arbitrary small set without
changing appropriate $x\in\QQ$. Allowing some inaccuracy we'll talk about hyperrational number $x=x_n$.
Let us call subsequence $x_{n_k}=x_J$ of $x_n$ \emph{big subsequence} if the 
subset $J=\{n_k: k\in\omega\}$ of indexes is big.

The hyperrational number $x=x_n$ is infinitely large  if and only if any big subsequence $x_{n_k}$ is unbounded,
is infinitesimal if and only if $1/x$ is infinitely large.

\section{Propositions}
\begin{theor} \label{tmono}
  For any $x=x_n\in\QQ$, there exists big subsequence $x_{n_k}$, which is strictly increasing or strictly decreasing or stationary.
  The cases are mutually exclusive. In case of $x<\infty$ the subsequence is fundamental and any two such subsequences $x_{n_k}$
  and $x_{m_k}$ are equivalent in traditional metric sense $\lim\limits_{k\to\infty}(x_{n_k}-x_{m_k})=0$
\end{theor}

Hyperrational numbers given by increasing (decreasing) sequences we call \emph{left} (\emph{right}) numbers.
Hyperrational numbers given by stationary sequences are exactly rational numbers.

Let $\QF=\{x\in\QQ~|~x<\infty\}$ be the set of all finite hyperrational numbers and infinitesimals.
\begin{theor}\label{treal}
  The set $\QF$ is local ring, whose unique maximal ideal is the set of infinitesimals $\II$. 
  The factor ring of $\QF$ modulo $\II$ is the field isomorphic to the field of real numbers:
  $$\QF / \II\simeq\R$$
\end{theor}

\section{Proofs}
We give the proof of the theorem \ref{t0} which is a bit more simple than the one from \cite{com_neg} 
in part of implication "selective $\Rightarrow$ normal".

\begin{proof}[Proof of theorem \ref{t0}]
  We use following proof schema: selective $\Rightarrow$ normal $\Rightarrow$ Ramsey $\Rightarrow$ selective.

  selective $\Rightarrow$ normal. Let $\{A_i\}_{i\in\omega}\subset\FF$ be arbitrary collection of big sets.
  We can assume with no loss of generality that $\forall i\in A_k~:~ i>k$, 
  because we can replace $A_k$ with big subsets $A_k'=\{i\in A_k | i>k\}$.
  Let us define a mapping  $f:\omega\to\omega, f(i)=\min\{j\in\omega~|~i\notin A_j\}$.
  Thus $f(i)\leq i$ because $i\notin A_i$. 
  Obviously, $f^{-1}(j)\cap A_j = \varnothing$ and, so, sets $f^{-1}(j)$ are small for all $j$.
  Then there exists big set $B$ such that $f$ is injective on $B$.

  We'll construct big set $A$ with the property: $\forall i,j\in A : i<j \Rightarrow i<f(j) \leq j$. 
  Such a set satisfies the conditions of the statement.

  Now let us construct subsets $P_k$ of $B$ as follows:
  \begin{eqnarray*}
    m_0&=&f(b_0)=\min f(B), P_0=\{b_0\},\\
    S_k &=& \{s\in f(B) ~|~ s>\max\bigcup\limits_{i=0}^{k-1}P_i\},~k\geq 1,\\
    m_k&=&f(b_k)=\min S_k,~k\geq1, \\
    P_k&=&\{b\in B ~|~ m_{k-1}< f(b) \leq m_k\},~k\geq 1
  \end{eqnarray*}

  All $P_k$ are finite because $f|_B$ is injective. Obviously, $\cup_{i=0}^k P_i=\{b\in B~|~f(b)\leq m_k\}$. 
  All $S_k$ are infinite because $f(B)$ is infinite. Thus $S_k\neq\varnothing$ and all $m_k$ are correctly defined.
  Note that $S_{k+1}\subseteq S_k$ and so $m_k\leq m_{k+1}$ for all $k$.
  
  In fact $m_k=f(b_k)\leq b_k<m_{k+1}$ because $b_k\in P_k$ and $m_{k+1}>\max \cup_{i=0}^k P_i$.
  So, the sequence $m_k$ is strictly increasing.
  The sets $P_k$ are not empty because at least $b_k\in P_k$.

  Now if $y\in P_{k+2}$ then $\max(\cup_{i=0}^k P_i)<m_{k+1}<f(y)$ and so $x<f(y) \leq y$ for all $x\in\cup_{i=0}^k P_i$.

  We have $\coprod_k P_k = B$. One of two sets $\coprod_k P_{2k}$ and $\coprod_k P_{2k+1}$ 
  is big and partitioned with small sets $P_n$ where n is odd or even. Let $A$ be big selection from this partition.
  Let $i<j$ be arbitrary elements of $A$. Then $i\in P_k$ and $j\in P_{k+2s}$ for some $k$ and $s>0$. So, $i<f(j)\leq j$
  and $j\in A_i$.

  normal $\Rightarrow$ Ramsey. Let $[\omega]^2=P\coprod Q$ be some partition of $[\omega]^2$.
  We consider following subsets of $\omega$:
  \begin{eqnarray*}
    P_i=\{j\in\omega~|~\{i,j\}\in P~\text{and}~j>i\} \\
    Q_i=\{j\in\omega~|~\{i,j\}\in Q~\text{and}~j>i\} 
  \end{eqnarray*}

  Obviously, $\omega=P_i\coprod Q_i\coprod \{1,\dots,i\}$ for all $i$. So, for fixed $i$ one and only one of $P_i$ and $Q_i$ is big.
  Let $B=\{i\in\omega~|~P_i\in\FF\}$ and $C=\{i\in\omega~|~Q_i\in\FF\}$. One of $B$ and $C$ is big. Let it be $B$. So, we have
  family $\{P_i\}_{i\in B}$ of big sets. By the definition of normal ultrafilter there exists big set $A$ such that for any
  two elements $i<j$ from $A$ we have $j\in P_i$ which means $\{i,j\}\in P$. So, $A$ is homogeneous.
  
  Ramsey $\Rightarrow$ selective. Let $\{S_i\}_{i\in\omega}$ be a small partition of $\omega$. Consider 
  $Q=\{\{i,j\}\in[\omega]^2~|~\exists k: i\in S_k~\text{and}~j\in S_k\}$ and $P = \omega - Q$. There exists
  big subset $H\subset\omega$ such that $[H]^2\subset Q$ or $[H]^2\subset P$. But $[H]^2\subset Q$ implies $H\subset S_k$
  for some $k$ which is impossible because $S_k$ is small. Thus, $[H]^2\subset P$ and the intersection $H\cap S_k$ can not
  contain more than one element for any $k$. We can add elements to $H$ if some of the intersections are empty. So, $H$ is
  the desired big selection.
\end{proof}

\begin{lemma}\label{l2}
  For any injection $\pi:\omega\to\omega$ there exists big subset $B$ such that $\pi|_B$ is increasing.
\end{lemma}
\begin{proof}
  We define the partition $[\omega]^2=P\coprod Q$ as follows 
  $P = \{\{i,j\}\in[\omega]^2~|~i<j~\text{and}~\pi(i)<\pi(j)\}$ and $Q=[\omega]^2 - P$.
  There exists homogeneous big set $B$ for this partition. But $[B]^2$ can not be subset of $Q$ because
  there is no infinitely decreasing sequences in $\omega$. So, $[B]^2\subset P$ and $\pi|_B$ is increasing.
\end{proof}

\begin{proof}[Proof of theorem \ref{tmono}]
  
  Let $x=x_n\in\QQ$ be arbitrary hyperrational number. 
  First of all let us consider the equivalence relation on $\omega:~n\sim k$ if and only if
  $x_n=x_k$. If the partition corresponding to this relation have big subset then there is big stationary subsequence
  of $x_n$.

  Otherwise, we consider $D\subset\omega$ be big selection from this partition. 
  So, for $k,n\in D$ if $k\neq n$ then $x_k\neq x_n$.
  For any $k\in\Z$ we define subsets $I_k = \{n\in D~|~x_n\in(k,k+1]\}$.
  It is clear that choosing nonempty subsets $I_k$ we  get the partition of $D$. 
  If this partition is small then there exists big selection $B$. One and only one of two sets
  $B_1=\{n\in B~|~ x_n>0\}$ and $B_2=\{n\in B~|~ x_n<0\}$ is big. 
  Note that the set $\{x_n~|~n\in B_i\}$ for big $B_i$ has order type $\omega$ in case of $B_1$
  and $\omega^*$ in case of $B_2$. 

  Otherwise, $I_k$ is big for some $k$.
  Thus, $x_{I_k}$ is big bounded subsequence and number $x$ is finite or infinitesimal.
  Only one of two sets $E_1=\{n\in B~|~ x_n<x\}$ and $E_2=\{n\in B~|~ x_n>x\}$ is big.
  
  Let show that if $E_1$ is big then there exists big subset $B_3\subset E_1$ such that
  $\{x_n~|~n\in B_3\}$ has order type $\omega$ and $x_{B_3}$ is fundamental.
  For this purpose let us define the partition of $E_1$
  as follows: $J_s = \{i\in E_1~|~x-\frac{1}{s}\leq x_i<x-\frac{1}{s+1}\}$.
  $J_s$ subsets are small for all $s$ because if $J_s$ is big for some $s$ then 
  the number $x'$ given by $x_{J_s}$ is equal to $x$ but on other hand $x'<x$. Contradiction.
  Thus, we can get big selection $B_3$ from the partition $\{J_s\}_{s<\omega}$.
  The sequence $x_{B_3}$ is obviously fundamental and the set $\{x_n~|~n\in B_3\}$ has order type $\omega$.

  Similarly If $E_2$ is big we can get big subset $B_4\subset E_2$ such that  
  the sequence $x_{B_4}$ is fundamental and the set $\{x_n~|~n\in B_4\}$ has order type $\omega^*$.
  
  If $x_K$ and $x_L$ are big fundamental subsequences then $x_{K\cap L}$ is big fundamental
  subsequence of both $x_K$ and $x_L$. Thus, $x_K$ and $x_L$ are equivalent in traditional metric sense.

  Let $C$ be the only big subset from subsets $B_i$ and $X=\{x_n\}_{n\in C}\subset\Q$ 
  subset of sequence elements. We define injection $\pi:\omega\to\omega$ as follows:
  $$
    \pi(n)=
      \begin{cases}
        \text{index of}~\min (X-\{x_{\pi(1)},\dots,x_{\pi(n-1)}\}),~\text{if}~C=B_1~\text{or}~C=B_3\\
        \text{index of}~\max (X-\{x_{\pi(1)},\dots,x_{\pi(n-1)}\}),~\text{if}~C=B_2~\text{or}~C=B_4\\
      \end{cases}
  $$
  
  Thus, $\pi$ is descending or ascending ordering of $X$ and for any $i<j$ we have $x_{\pi(i)}<x_{\pi(j)}$ in first case
  and $x_{\pi(i)}>x_{\pi(j)}$ in second case. According to lemma \ref{l2} there exists big subset $E\subset\omega$
  such that $\pi_E$ is increasing and subsequence $x_{E\cap C}$ is monotonic.
\end{proof}

We call $\nu(x)\in\R$ the \emph{value} of number $x\in\QF$. 

\begin{proof}[Proof of theorem \ref{treal}]
  Let us define mapping $\nu:\QF\to\R$. For any $x=x_n\in\QF$ we set $\nu(x)$ equal to limit of some big fundamental subsequence
  of $x_n$. This limit is uniquely defined as follows from the theorem \ref{tmono}. It is easy to see that $\nu$ is epimorphism of
  $\Q$-algebras and $\ker\nu=\II$ is the ideal of all infinitesimals in $\QF$. All elements of compliment of $\II$ in $\QF$
  are invertible. Thus, $\II$ is only maximal ideal in local ring $\QF$ and $\QF/\II\simeq\R$.
\end{proof}


\begin{thebibliography} {99}
  \bibitem{rob0}
    Robinson, Abraham \emph{``Non-standard analysis''}, Proc. Roy. Acad. Sci. Amst., ser. A 64 (1961), 432-440
  \bibitem{rob1}
    Robinson, Abraham \emph{``Non-standard analysis''}, Princeton University Press, 1996, ISBN: 978-0-691-04490-3
  \bibitem{jech}
    Jech, Thomas \emph{``Set theory''}, Springer, 2006, ISBN 3-540-44085-2
  \bibitem{shelah}
    Wimmers, Edward (March 1982), \emph{``The Shelah P-point independence theorem''}, 
    Israel Journal of Mathematics, Hebrew University Magnes Press, 43 (1): 28–48, doi:10.1007/BF02761683
  \bibitem{com_neg}
    Comfort W. W., Negrepontis S., \emph{``The theory of ultrafilters''}, 
    Springer-Verlag Berlin Heidelberg, 1974, ISBN-10: 0387066047
\end{thebibliography}
\end{document}